\documentclass[12pt]{article}

\usepackage{amsmath,amssymb}
\usepackage{graphicx}

\begin{document}

\begin{center}
{\bf Bayes factors with (overly) informative priors}\\[18pt]

Richard Lockhart\footnote{The author thanks  Michael 
Stephens for many useful conversations on the topics discussed here.
The author acknowledges grant support from the  
Natural Sciences and Engineering Research Council of Canada.} \\[12pt]
Department of Statistics \& Actuarial Science\\
Simon Fraser University 
\\
Burnaby, BC \ V5A 1S6 CANADA
\end{center}

\centerline{\bf Abstract}

{\it Priors in which a large number of parameters are specified to be independent are dangerous; they make it hard to learn from data.  I present a couple of examples
from the literature and work through a bit of large sample theory to show what happens.}

\bigskip

\noindent{{\sl Keywords}:  Teaching examples, hierarchical priors, independence priors, pitfalls in prior specification.

\medskip

{\sl MSC Classification Code}\/: 62F15, 62E20}

\section{Introduction}

B\'{e}lisle et al.~(2002) fit five interesting models to a data set on progression of Alzheimer's disease.
They use Bayes factors (both prior and posterior) for model selection.  Most of the Bayes factors reported are 
numbers of the form $10^x$ with $x$ on the order of plus or minus tens to hundreds; that is, in most cases
when two models are compared one of the two is overwhelmingly preferred.  This is not too surprising from a Bayes
perspective but the results reveal some comparisons which are surprising for a frequentist.

In some of these models there is a parameter which varies
from individual to individual.  In other models this parameter is constant across individuals.  For a frequentist
the latter model is a submodel of the former.  A frequentist would,  I submit, be surprised to be told that in the
face of the data the submodel is overwhelmingly (like $10^{12}$ time more likely) to be preferred to the richer
model.  In the results presented by B\'{e}lisle et al.~(2002) this is exactly what happens: their model 5 
is strongly preferred to their model 4 of which model 5 is a submodel.  They report prior Bayes factors of $10^{22.7}$ and $10^{12}$ in favour of the submodel in two different data sets. They also compute  
posterior Bayes factors, known to be computationally easier,  of $10^{26.5}$ and $10^{16.3}$ in favour of the submodel in the same two datasets.

Their modelling suggests, in fact, that neither of these models is as good as their models 1, 2 or 3 which are more
complex.  It seems nevertheless worthwhile to understand the potential for Bayes factors strongly to prefer
submodels to full models.  In this note I use the one way layout model to give an example in which a wrong submodel
is overwhelmingly preferred to the full model.  

Priors in which a large number of parameters are specified as being independent can easily cause problems. For a second example I consider a stylized survey sampling problem rom Wasserman (2004) who presents it as a problem in which Bayes methods struggle.

In Section~\ref{sec:ksample} I  present the mathematical details of the one way layout problem and discuss briefly the relation to the situation in B\'{e}lisle et al.~(2002).  In Section~\ref{sec:wasserman} I simplify Wasserman's problem and present the standard Bayesian analysis.  I finish the paper with a short discussion arguing that

\begin{enumerate}

\item Priors in which many parameters are independent can be too informative to be safely used in data analysis.

\item Hierarchical priors, as Bayesians know well, avoid these pitfalls.

\item In order to learn from what happens in one measurement about another measurement a Bayesian must, before  making the first measurement, regard the two outcomes as dependent.

\end{enumerate}

I hope the examples here might be useful for pedagogy and for highlighting dangers of careless use of priors.

\section{Example 1: the one way layout}\label{sec:ksample}

Consider the $k$ sample problem with known variance.  We suppose we have data $X_{ij};j=1,\ldots,n_i;i=1,\ldots,k$.
We will be considering the standard analysis of variance model.  We assume  
the $X_{ij}$ to be independent with $X_{ij}$ having a normal distribution with mean $\mu_i$ and variance 1.
This will be Model 1.  Model 2 is the nested submodel with all $\mu_i=0$; that is, for Model 2 the data are iid
standard normal.  

This is usually treated by frequentists as a hypothesis testing problem but we will look at it here as a
model selection problem.  We will consider a Bayesian approach.  Conditional on Model 2 holding there are no
further parameters on which to put a prior distribution.  Conditional on  Model 1 holding
we need a prior for the $\mu_i$.  We consider the simple conjugate prior under which the means $\mu_i$ are 
independent $N(0,\tau^2)$ random variables.  

If the 2 models are given prior probabilities $\pi_1$ and $\pi_2=1-\pi_1$ then the posterior probability that model
2 holds given the data is
$$
\frac{\pi_2}{\pi_1 + \pi_2 F}
$$
where $F$ is the prior Bayes factor given by
$$
F = \frac{f_2(x)}{f_1(x)}
$$
with $f_i$ being the marginal density of the data under model $i$. (That is, $f_1$ is the model joint density averaged with
respect to the prior on the $\mu_i$ and $f_2$ is the joint density of iid standard normals.)

Elementary normal distribution theory calculations show that 
$$
f_1(x) = (2\pi)^{-n/2} \prod_{i=1}^k (1+n_i \tau^2)^{-1/2} 
\exp\left\{-\frac12\sum_{ij} x_{ij}^2 + \frac{\tau^2 \sum_i \left(\sum_j x_{ij}\right)^2}{2(1+n_i\tau^2)}\right\}
$$
where we let $x$ be the vector of all the $x_{ij}$ values.
It follows that
$$
\log(F) = -\sum_i \frac{ \tau^2 \left(\sum_j X_{ij}\right)^2}{2(1+n_i \tau^2)} +\frac{1}{2} \sum_i \log(1+n_i\tau^2).
$$
Now we want to consider the frequentist properties of $\log(F)$ by imagining a sequence $\mu_1, \mu_2, \ldots$. I will take these values to be deterministic 
and try to see how $\log(F)$ behaves when there are many samples, that is, when $k$ is large. I begin by computing the first few moments of $\log(F)$.  For simplicity I will take all the $n_i$ to be equal and use $n$ to denote their common value.  In this case, for each $i$ we find that 
$$
T_i \equiv \sum_j X_{ij} \sim \text{Normal}(n\mu_i,n)
$$
and these sums, $T_i$, are independent over $i$.  Thus
$$
{\rm E}\left(T_i^2\right) = n+n^2 \mu_i^2.
$$
Thus the mean of $\log(F)$ is 
$$
{\rm E}\left[\log(F)\right] = -\frac{ \tau^2}{2(1+n\tau^2)} \left( nk +n^2 \sum_i \mu_i^2\right)  +\frac{1}{2} k \log(1+n\tau^2).
$$

We can now check that if $\overline{\mu^2} = \sum_i \mu_i^2/k$ is small enough then $\log(F)$ can easily be positive. A simpler computation ensues if the $\mu_i$ are in fact generated in an iid way from a Normal$(0,\epsilon^2)$ distribution; I now make this assumption. Then holding $n$ and $\tau$ fixed and letting $k$, the number of samples, tend to $\infty$, we see that 
$$
2\frac{\log(F)}{k} \to \log(1+a) - \frac{a(1+n\epsilon^2)}{1+a} 
$$
where $a=n\tau^2$. Define $f(a)=\log(1+a)-a/(1+a) $ and check that $f(0)=0$ and $f'(a) = a/(1+a)^2>0$. The term $a n\epsilon^2/(1+a)$ is no more than $n\epsilon^2$ which converges to 0 as $\epsilon\to 0$. So for each $n,\tau$ pair there is an $\epsilon$ so small that the limit of $2\log(F)/k$ is positive.  For this sequence --- a sequence for which the null hypothesis is clearly false --- the Bayes factor in favour of the null model, grows exponentially fast.  If $n=10$ and $\tau = 1$, for instance, then $a=10$.  Numerically we find that with $\epsilon=0.3$ the limit of $\log(F)/k$ is $0.67\approx 2/3$.  In other words if the actual variability in the means is 30\% of the variability predicted by the prior then the submodel is roughly $\exp(2k/3)$ times as credible as the full model.

\begin{figure}\label{MedianPlot}
\includegraphics[width=\textwidth]{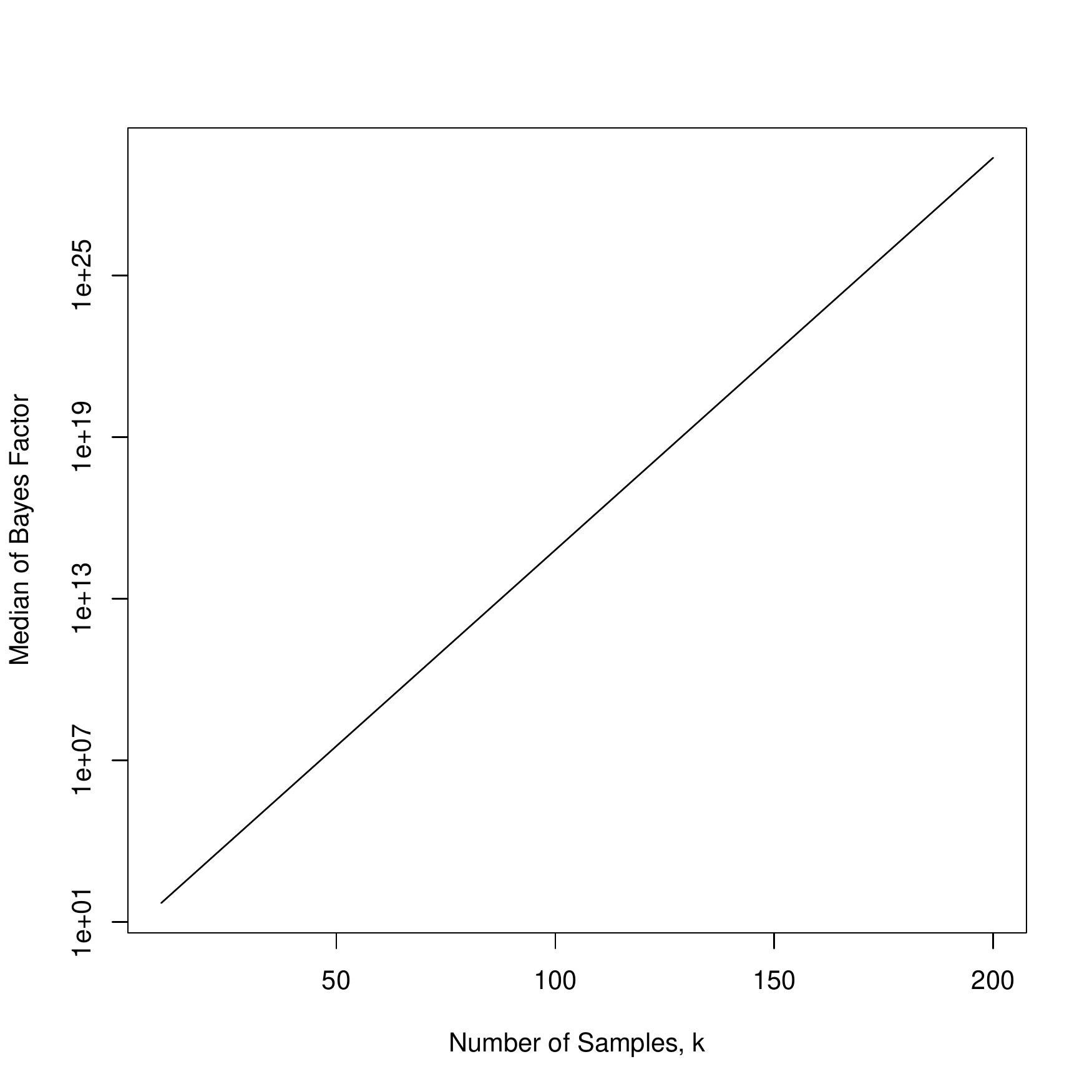}
\caption{The median of $F$ is plotted against the number $k$ of samples for $n=10$, $\tau=1$ and $\epsilon=0.3$. Note the logarithmic scale on the vertical axis.}
\end{figure}

A reasonable form of the analytic law of $\log(F)$ is available when the $\mu_i$ are iid $N(0,\epsilon^2)$. In this case
$$
T_i = n\mu_i + \sum_j (X_{ij}-\mu_i)
$$
has a normal distribution with mean 0 and variance 
$$
\sigma^2 = n^2 \epsilon^2 + n.
$$
Since the $T_i$ are iid it follows that 
$$
\sum_i T_i^2/\sigma^2 \sim \chi^2_k.
$$
Thus 
$$
\log(F) = k\log(1+n\tau^2)/2 - \tau^2 \sigma^2 \sum_i T_i^2/\sigma^2
$$
so that 
$$
P(\log(F) > kt) = P\left[\frac{\chi_k^2}{k} > (1+n\tau^2)\frac{2t-\log(1+n\tau^2)}{\sigma^2 \tau^2}\right]
$$
Using this form I plot, in Figure~\ref{MedianPlot},  the median of $F$ against the number $k$ of samples for $n=10$, $\tau=1$ and $\epsilon=0.3$.  I use a logarithmic scale for $F$ on the $y$-axis and the nearly perfect linearity is evident. Note, too, that the actual median values are very large by the time $k=50$ and vast for $k=200$. 

The situation here is relatively simple.  In the model above the value of $\mu_i$ is estimated from 10 observations. The prior specified iid $\mu_i$ with variance 1 so that the standard deviation of the $\mu_i$ arising was \emph{a priori} expected to be 1. In fact that standard deviation was smaller, 0.3 only.  As a consequence when $k$ is moderately large the prior for Model 1, the less restrictive model, is contradicted for the data.  Model 2 is also contradicted (and false).  For $\epsilon<0.404$ or so, as $k\to\infty$ with $n$ and $tau$ held fixed, we find that $F$ tends exponentially fast to $\infty$ and we erroneously conclude that the means are all 0.

In B\'{e}lisle et al.~(2002), Model 5 differs from Model 4 in allowing a certain variance parameter to be different from subject to subject.  These variance parameters were taken to be independent across individuals and to have expectation 2; they were given inverse Gaussian distributions with this expectation and with standard deviation 2.  There were $n=65$ subjects in one of the two data sets and $n=55$ in the other.   I do not have access to the data (honestly, I haven't tried to get access to the data) so I cannot check, but I assume that it is likely that a fit with some smaller mean and or some other standard deviation for these variance parameters would yield more natural Bayes factors.

The phenonemon described above arises in other model selection problems.  For instance, we regress $Y$ on $p$ predictors $X_1,\ldots,X_p$ and assume that each $X_i$ has a non-zero coefficient independently of all the other covariates and that $\theta \in(0,1)$ is the common value then if we do model selection by Bayesian methods the posterior will be concentrated on the event that the number of active predictors is in the range $p\theta \pm 2 \sqrt{p\theta(1-\theta)}$.

\section{Example 2: Survey Sampling}\label{sec:wasserman}

A related phenomenon occurs in Larry Wasserman's book \emph{All of Statistics}.  Example 11.19  in the book is a description of a missing data problem but it has the flavour of survey sampling. There are $B$ parameters $\theta_1,\ldots,\theta_B$, each a probability in $[0,1]$. We draw a sample with replacement from the index set $\{1,\cdots,B\}$; let $X_1,\ldots,X_n$ be the indices selected. It is assumed that these are iid and uniform on $\{1,\cdots,B\}$.  Having selected, on the $i$th draw, unit $X_i$, we may or may not observe an observation $Y_i$ which is Bernoulli with parameter $\theta_{X_i}$.  The example given supposes that `may or may not' is determined by a Bernoulli variable $R_i$ with parameter $\xi_{X_i}$ which is known to the analyst.  In the text it is argued that the Horwitz-Thompson estimator of 
$$
\psi_B \equiv \frac{\sum_1^B \theta_j}{B}
$$ 
is a good one and that Bayesian methods fail.  Here I simplify the model, eliminating the possibility of missing data and make the sampling law a bit more in line with survey sampling work.

I will assume that we generate a set $J$, the sample, of indices in $\{1,\ldots,B\}$ with design probability $\pi_J$; I assume the design is not informative; that is the $\pi_J$ are known to the surveyor in advance and unrelated to $\psi$. Then for each $j\in J$ we observe $Y_j$ with a Bernoulli$(\theta_j)$ distribution independently of all other $Y$ values.  I let $Y$ denote the vector of $J$
observed values.

The likelihood is the probability of observing the data we observed:
$$
\pi_J \prod_{j \in J}\theta_j^{Y_j}(1-\theta_j)^{1-Y_i}.
$$
This  likelihood is observed to depend  on only  $n$  of the $B$ entries in the vector $\theta$ of all $\theta_i$.  Wasserman  now argues that ``for most $\theta_j$ the posterior distribution is equal to the prior distribution since those $\theta_j$ do not appear in the likelihood.''   The problem, however, is that the prior is not clearly specified.  If I specify that the $\theta_i$ are \emph{a priori} iid with density $\pi(\cdot)$ then the posterior for unobserved $\theta_i$ is still iid with density $\pi(\cdot)$. This appears to me to be the argument intended by Wasserman.  But this prior specifies that if I measure one $\theta_i$ I do not learn about any other $\theta_i$. This sort of prior is exactly the sort discussed above and has the same weakness.  The unobserved $B-n$ parameters $\theta_j$ will average,
essentially, to $\psi\equiv \int_0^1 \theta \pi(\theta) d\theta$.  If, as the text assumes, $B-n$ is very large compared to $n$ then of course I will not learn anything important about $\psi_B$; \emph{a priori} I knew that $\psi_B$ was \emph{close} to $\psi$.

I would argue that the correct message here is as above; independence priors about many parameters are priors which \emph{deny learning}. If I had to guess the mass of a flea I would struggle  but if you let me weigh one of a group of many I would suddenly know far more about the average weight of the group of fleas.  I need a hierarchical prior.

In the Wasserman example (as modified here) one might specify such a hierarchical prior by taking the $\theta_i$ to be iid with a Beta($\alpha,\beta$) distribution and then putting  a prior on $(\alpha,\beta)\in(0,\infty)^2$.  For the discussion which follows I reparametrize the Beta distribution in terms of its mean
$$
\psi=\frac{\alpha}{\alpha+\beta}
$$
and its variance 
$$
\eta\equiv \frac{\psi(1-\psi)}{\alpha+\beta+1}.
$$
The parameters $\alpha$ and $\beta$ may be expressed in terms of $\psi$ and $\eta$ by noting
$$
\alpha+\beta = \frac{\psi(1-\psi)}{\eta} - 1
$$
and
$$
\alpha = (\alpha+\beta) \psi = \frac{\psi^2(1-\psi)}{\eta} -\psi
$$
This gives
$$
\beta =  \frac{\psi(1-\psi)}{\eta} - 1 -( \frac{\psi^2(1-\psi)}{\eta} -\psi)= \frac{(1-\psi)(\psi-\psi^2-\eta)}{\eta} .
$$

A simple prior specification  might then be to give $\psi$  some Beta($\alpha_0,\beta_0$) density denoted $h_0(\psi)$ on $[0,1]$ and to give $\eta$ some conditional density given $\psi$.  Notice that $0 \le \eta \le \psi(1-\psi)$ so this conditional density must be supported on this interval. One simple conditional prior would be to take $\eta$ to be uniform on $[0,\psi_\infty(1-\psi_\infty)]$; I use $h(\eta|\psi)$ for some generic conditional prior for $\eta$ given $\psi$ and then specialize where needed to the uniform case where
$$
h(\eta|\psi) =\frac{1\left(0 \le \eta \le \psi(1-\psi)\right)}{ \psi(1-\psi)}.
$$
 \emph{A priori} we have ${\rm E}\left(\psi_B|\psi,\eta \right)=\psi$.

The joint law of the data, the parameters $\theta_i$ and the hyperparameters $\psi$ and $\eta$ is 
\begin{equation}\label{eq:jointall}
\pi_J \prod_{j \in J} \theta_j^{Y_j}(1-\theta_j)^{1-Y_j} \prod_j \frac{\theta_j^{\alpha-1}(1-\theta_j)^{\beta-1} }{{\rm Beta}(\alpha,\beta)}
h(\eta|\psi)h_0(\psi).
\end{equation}
From this I now deduce: the joint law of the data $J,Y$ and the hyperparameters $\psi$ and $\eta$; the usual Bayes estimate of $\psi$, namely ${\rm E}(\psi|Y)$; and the Bayes estimate of $\psi_B$, namely ${\rm E}(\psi_B|Y)$.

For the first of these I simply integrate away all the $\theta_j$.  For $j\not\in J$ we have
$$
\int_0^1  \frac{\theta_j^{\alpha-1}(1-\theta_j)^{\beta-1} }{{\rm Beta}(\alpha,\beta)} d\theta_j = 1.
$$
For $j \in J$ the integral needed becomes
$$
\int_0^1 \theta_j^{Y_j}(1-\theta_j)^{1-Y_j}   \frac{\theta_j^{\alpha-1}(1-\theta_j)^{\beta-1} }{{\rm Beta}(\alpha,\beta)} d\theta_j 
= \frac{{\rm Beta}(\alpha+Y_j,\beta+1-Y_j)}{{\rm Beta}(\alpha,\beta)}.
$$
Considering the two cases $Y_j=1$ and $Y_j=0$ separately we find
$$
\int_0^1 \theta_j^{Y_j}(1-\theta_j)^{1-Y_j}   \frac{\theta_j^{\alpha-1}(1-\theta_j)^{\beta-1} }{{\rm Beta}(\alpha,\beta)} d\theta_j  = \psi^{Y_j} (1-\psi)^{1-Y_j}.
$$
Thus the joint law of $J,Y,\psi,\eta$ is
$$
\pi_J\prod_{j\in J} \psi^Y_j (1-\psi)^{1-Y_j} h(\eta|\psi)= \pi_J \psi^S(1-\psi)^{|J|-S} h(\eta|\psi)h_0(\psi)
$$
where $S=\sum_J Y_j$ and we use $|J|$ to denote the cardinality of $J$.

The conditional density, $f(\psi,\eta|Y,J)$, of $\psi,\eta$ given $Y$ and $J$ is then
$$
\frac{\psi^S(1-\psi)^{|J|-S} h(\eta|\psi)}{{\rm Beta}(S+1,|J|-S+1)}.
$$
Now integrate out $\eta$ to find the conditional density of $\psi$ given $Y$ and $J$, namely,
$$
\frac{\psi^S(1-\psi)^{|J|-S}}{B(S+1,|J|-S+1)}h_0(\psi).
$$
If $h_0$ is the Beta($\alpha_0,\beta_0$) then this posterior is a Beta density; the Bayes estimate of $\psi$ is the mean of this density, namely,
$$
\hat\psi = \frac{S+\alpha_0}{|J|+\alpha_0+\beta_0}
$$
I remark that the improper prior $\psi^{-1}(1-\psi)^{-1}$ arises in the limit $\alpha_0\to 0$ and $\beta_0\to 0$.  This  leads to the Horwitz-Thompson estimator
$$
\hat\psi_{\rm HT} = \frac{S}{|J|}
$$
suggested by Wasserman.  Of course this is our estimate of $\psi$ while Wasserman is estimating $\psi_B$.

Finally I compute the Bayes estimate of $\psi_B$ which is
$$
{\rm E}(\psi_B|Y,J)  = \frac{1}{B} \sum_j {\rm E}(\theta_j|Y,J).
$$
I compute ${\rm E}(\theta_j |Y)$ separately according to $j\not\in J$ or $j\in J$.

Notice that
$$
{\rm E}(\theta_j|Y,J) = {\rm E} \left\{{ \rm E}\left( \theta_j |  Y, J,\psi,\eta \right) | Y,J \right\} .
$$
In order to compute the inside conditional expectation I return to joint law~(\ref{eq:jointall}).  For $j\not\in J$ this joint law may be written in the form
$$
\frac{\theta_j^{\alpha-1}(1-\theta_j)^{\beta-1} }{{\rm Beta}(\alpha,\beta)} \times C
$$
where the quantity $C$ does not depend on $\theta_j$.  It follows that given $Y,J,\psi,\eta$ the conditional law of $\theta_j$ is, for $j\not\in J$, Beta with parameters $\psi$ and $\eta$. Thus 
$$
{\rm E}(\theta_j |Y,J,\psi,\eta) =\psi .
$$
Finally we see that 
$$
{\rm E}(\theta_j|Y,J) = {\rm E}(\psi |Y,J) = \hat\psi.
$$

When $j\in J$ the kernel of~(\ref{eq:jointall}) is
$$
\theta_j^{\alpha+Y_j-1}(1-\theta_j)^{\beta+1-Y_j+1}
$$
which is the kernel of the Beta density with parameters $\alpha+Y_j$ and $\beta+1-Y_j$.  The mean of this density is 
$$
\frac{\alpha+Y_j}{\alpha+\beta+1}
$$
and we need to compute
$$
{\rm E}\left\{{\rm E}\left(\theta_j| Y, J,\psi,\eta\right) | Y,J \right\} 
$$
This is
$$
\int {\rm E}\left(\theta_j| Y, J,\psi,\eta\right) h(\eta|\psi)f(\psi|Y,J) d\eta d\psi.
$$
The required conditional density is given above but we must write the inner conditional mean in terms of $\psi$ and $\eta$. First
we get
$$
\frac{\alpha}{\alpha+\beta+1} =   \psi(\alpha+\beta)\frac{\eta}{\psi(1-\psi)}=  \psi\left\{ \frac{\psi(1-\psi)}{\eta} - 1\right\}
\frac{\eta}{\psi(1-\psi)}=\psi-\frac{\eta}{1-\psi}.
$$
We also find 
$$
\frac{Y_j}{\alpha+\beta+1}=Y_j \frac{\eta}{\psi(1-\psi)}.
$$
Our integral becomes
$$
\int \int  \left(\psi-\frac{\eta}{1-\psi} +Y_j\frac{\eta}{\psi(1-\psi)}\right) f(\psi,\eta|Y,J) d\eta d\psi.
$$
I now use the conditional uniform distribution of $\eta$ to do the inside integral to get
$$
\int \left( \psi - \frac{\psi}{2}+ \frac{Y_j}{2}\right) f(\psi|Y,J) d\psi = \hat\psi+ \frac{Y_j-\hat\psi}{2}.
$$

In summary we have found that
$$
{\rm E}(\theta_j |  Y,J) = \begin{cases} \hat\psi & j \not\in J \\ \hat\psi+\frac{Y_j-\hat\psi}{2} & j \in J \end{cases}.
$$
From this we find that the Bayes estimate of $\psi$ is 
$$
\hat\psi - \frac{|J|\hat\psi+S}{2B}  = \hat\psi + \frac{|J|}{B} \frac{\hat\psi_{\rm HT}-\hat\psi}{2}.
$$
When $|J|$ is small compared to $B$ the second term is negligible.  In fact, for the special choice of the non-informative improper prior $\psi^{-1}(1-\psi)^{-1}$ the second term vanishes exactly and the Bayes estimate is just the Horvitz-Thompson estimator.  Finally remark that the second term is 
$$
 \frac{|J|}{2B} \left\{\frac{S+\alpha_0}{|J|+\alpha_0+\beta_0} -\frac{S}{|J|}\right\} = \frac{(S+\alpha_0)|J|-S(|J|+\alpha_0+\beta_0)}{2B(|J|+\alpha_0+\beta_0)}.
 $$
This quantity is easily seen to be bounded, in absolute value, by 
$$
\frac{\alpha_0/2}{B}+\frac{(\alpha_0+\beta_0)/2}{B} =O(1/B).
$$
Thus the correction to $\hat\psi$  is negligible if $B$ is large.

As a practical matter this Horvitz-Thompson has variance which is acceptably small. Treating the $\theta_j$ as fixed parameters the variance is 
$$
\sum_J \theta_i(1-\theta_j)/|J|^2.
$$
I do not think this quantity admits a useful estimate.
\section{Discussion}

When many parameters are specified as \emph{a priori} independent, the prior is very informative about averages of functions of those parameters.  The result is that estimates can depend little on the data and model selection methods can have very poor frequency theory properties.  Analysts using Bayesian methods must be careful of such priors and ought, in many cases, to make the prior hierarchical, inducing dependence between the parameters and permitting learning about some parameters by measuring others.

\bigskip
\noindent
{\sc References}

\begin{description}

\item B\'elisle, P.; Joseph, L.; Wolfson, D.; and Zhou, X. (2002). Bayesian estimation of cognitive decline in patients with Alzheimer's disease.
 \emph{The Canadian Journal of Statistics}, {\bf 30}, 37--54.

\item Wasserman, Larry (2004).
{\em All of Statistics}, Springer: New York.
\end{description}

\end{document}